\documentclass{elsart}

\journal{Journal of Number Theory}

\usepackage{amssymb}
\usepackage{amsxtra}
\usepackage[mathscr]{eucal}
\usepackage{graphics}
\usepackage{epsfig}

\newtheorem{theorem}{Theorem}[section]
\newtheorem{lemma}[theorem]{Lemma}

\theoremstyle{definition}

\newtheorem{definition}[theorem]{Definition}

\newtheorem{example}[theorem]{Example}

\newenvironment{proof}{\smallskip\noindent{\sc{Proof.}}}{\hfill$\square$\medskip}

\def\z{\mathsf Z}
\def\ze{\mathsf Z_{\mathsf e}}
\def\zo{\mathsf Z_{\mathsf o}}
\def\bigo{\operatorname{O}}    
\def\ord{\operatorname{ord}}
\def\divides{{\mathchoice{\mathrel{\bigm|}}{\mathrel{\bigm|}}{\mathrel{|}}%
       {\mathrel{|}}}}

\def\notdivides{\mathrel{\kern-3pt\not\!\kern3.5pt\bigm|}}

\def\proofof#1{\smallskip\noindent\emph{Proof of #1}.}
\def\endproofof{\hfill$\square$\medskip}

\begin{document}

\begin{frontmatter}

\title{Primitive divisors of
elliptic divisibility sequences}

\author{Graham Everest\thanksref{bennet}},
\ead{g.everest@uea.ac.uk}
\author{Gerard Mclaren},
\author{Thomas Ward\corauthref{cor}}
\ead{t.ward@uea.ac.uk}
\corauth[cor]{Corresponding author.
School of Mathematics, University of East Anglia,
Norwich NR4 7TJ, UK.}
\thanks[bennet]{Our thanks go to Mike Bennett for help in applying his results,
and to anonymous referees for several improvements.}

\begin{abstract}{Silverman proved the analogue of Zsigmondy's
Theorem for elliptic divisibility sequences. For elliptic curves in
global minimal form, it seems likely this result is true in a
uniform manner. We present such a result for certain infinite
families of curves and points. Our methods allow the first explicit
examples of the elliptic Zsigmondy Theorem to be exhibited. As an
application, we show that every term beyond the fourth
of the Somos-4 sequence has a primitive divisor.}
\end{abstract}

\end{frontmatter}

\section{Introduction}\label{intro}

Let~$A=\left(A_n\right)_{n\ge1}$ be an integer sequence.
A prime~$p$ dividing a term~$A_n$ is called a
\emph{primitive divisor} of~$A_n$ if~$p$ does not
divide any term~$A_m$,~$1\le m<n$. Thus, in the
list of prime factors of the terms of the sequence,
a primitive divisor is a new prime factor.
Sequences with the property that all terms (or all
terms beyond some point) have a primitive divisor
are of great interest.

\begin{definition}\label{def:z}
Let~$A=\left(A_n\right)_{n\ge1}$ be an integer sequence. Define
$$
\z(A)=
\max\{n\mid A_n\mbox{ does not have a primitive divisor}\}
$$
if this set is finite, and~$\z(A)=\infty$ if not. The
number~$\z(A)$ will be called the \emph{Zsigmondy bound} for~$A$.
\end{definition}

A striking early result is that of Zsigmondy~\cite{zsigmondy}.
For the Mersenne sequence
$$
M=\left(2^n-1\right)_{n\ge1},
$$
he showed
that
$$
\z(M)=6.
$$
More generally, Zsigmondy also showed that for any coprime
integers~$a$ and~$b$,
$$
\z\left(\left(a^n-b^n\right)_{n\ge1}\right)\le6.
$$
This line of development culminated in a deep result due to Bilu,
Hanrot and Voutier~\cite{MR2002j:11027}: for any non-trivial
Lucas or Lehmer
sequence~$L$,
$$
\z(L)\le30.
$$
Much of the arithmetic of linear recurrence
sequences extends to elliptic and bilinear
recurrence sequences (see~\cite[Chap.~10]{MR2004c:11015}
for an overview), and it is natural to ask if results
like that of Zsigmondy might hold for elliptic divisibility
sequences.

Let~$E$ denote an elliptic curve
defined over~$\mathbb Q$, given
in generalized Weierstrass form, and
suppose~$P=(x(P),y(P))$ denotes
a non-torsion rational point on~$E$
(see~\cite{MR92k:11058},~\cite{arithmetic},~\cite{MR87g:11070}
or~\cite{MR93g:11003}
for background on elliptic curves). For any
non-zero~$n\in \mathbb Z$, write
\begin{equation*}
x(nP)=\frac{A_n}{B_n},
\end{equation*}
in lowest terms, with~$A_n\in\mathbb Z$
and~$B_n\in\mathbb N$.
The sequence~$B_{E,P}=\left(B_n\right)_{n\ge1}$ is a divisibility
sequence, meaning that
$$
m\divides n\implies B_m\divides B_n.
$$
Such sequences have become
known as \emph{elliptic divisibility
sequences} (this terminology follows
a suggestion of Silverman; the term has also been used
for more general sequences related to rational
points on elliptic curves).
Silverman~\cite{MR89m:11027} showed that~$B_{E,P}$
satisfies an analogue of Zsigmondy's theorem.

\begin{theorem}{\sc [Silverman]}\label{ellipticzsig}
With~$E$ and~$P$ as above,
$$
\z(B_{E,P})<\infty.
$$
\end{theorem}

Our purpose here is to show that uniform explicit
bounds in Theorem~\ref{ellipticzsig} can be found
for certain infinite families of curves, after the
manner of~\cite{MR2002j:11027}. The methods allow
explicit versions of the theorem for particular
examples. Many of the bounds arrived at below can be
improved, similar methods may be applied to other elliptic
surfaces, and the
techniques used here may be applied to bound the
number of terms in an elliptic divisibility sequence
which are prime squares; further details in these directions may be found
in the thesis of the second named author~\cite{gerrythesis}.

\section{Main results}

The behaviour along the odd and even subsequences
of an elliptic divisibility sequence requires
slightly different treatment, so the following
refinement of Definition~\ref{def:z} will be
useful.

\begin{definition}\label{def:zoandze}
Let~$A=\left(A_n\right)_{n\ge1}$ be an integer sequence. Define
the \emph{even Zsigmondy bound}
$$
\ze(A)=
\max\{2n\mid A_{2n}\mbox{ does not have a primitive divisor}\}
$$
if this set is finite, and~$\ze(A)=\infty$ if not.
Similarly define the \emph{odd Zsigmondy bound}
$$
\zo(A)=
\max\{2n-1\mid A_{2n-1}\mbox{ does not have a primitive divisor}\}
$$
if this set is finite, and~$\zo(A)=\infty$ if not.
\end{definition}

Clearly~$\z(A)=\max\{\ze(A),\zo(A)\}$; in certain cases our
methods can bound explicitly either one of~$\ze$ and~$\zo$ but not
both.

\begin{theorem}\label{thetheorem}
Suppose the curve~$E$ is given by a
Weierstrass equation
\begin{equation*}
E:\quad y^2=x^3-T^2x,
\end{equation*}
with~$T>0$
square-free, and suppose that~$E$ has a
non-torsion point~$P$ in~$E(\mathbb Q)$.
Then
$$
\ze\left(B_{E,P}\right)\le10.
$$
If~$x(P)<0$, then
$$
\zo\left(B_{E,P}\right)\le3.
$$
If~$x(P)$ is a square,
then
$$
\zo\left(B_{E,P}\right)\le21.
$$
\end{theorem}

Notice that the existence of the point~$P$ certainly implies
that~$T\ge5$, so~$\log T$ is at least~$1.609$. This will be used several times
in the calculations below.

\begin{example}\label{y2=x3-25x}
Consider the curve
$$
E:\quad y^2=x^3-25x,
$$
with~$P=(-4,6)$.
We will show below that~$\z\left(B_{E,P}\right)=1$.
\end{example}

The assumption about~$T$ being square-free guarantees that~$E$ is
in global minimal form. Clearly an assumption of this kind is
necessary. It is always possible to clear arbitrarily many
denominators of the multiples~$x(nP)$ by applying suitable
isomorphisms, making
an explicit bound impossible.
Assuming the curve is in minimal form prevents this
possibility.

The most general form of result we can exhibit with our current
techniques will now be stated. Lang's Conjecture says that if~$E$
denotes an elliptic curve defined over~$\mathbb Q$ defined by a
Weierstrass equation in minimal form and if~$P$ denotes a
non-torsion rational
point on~$E$, then
\begin{equation}\label{langconjecture}
\hat h(P)\geq c\log \Delta(E).
\end{equation}
In~\eqref{langconjecture},~$\Delta(E)$ denotes the discriminant
of~$E$ and the constant~$c>0$ is uniform, independent of~$E$ and~$P$.
The family of curves in Theorem~\ref{thetheorem} is one for which
Lang's Conjecture is known to hold.

\begin{theorem}\label{thetheorem2}
Let~$\mathfrak F$ denote a family of elliptic curves~$E$, given by
Weierstrass models in global minimal form, and rational
points~$P,Q\in E(\mathbb Q)$, with~$P$ a non-torsion point and~$Q$
a~$2$-torsion point. Suppose that Lang's Conjecture holds for the
family; in other words, there is a uniform constant~$c=c(\mathfrak
F)>0$ such that for every triple~$(E,P,Q)\in \mathfrak F$, the
inequality~\eqref{langconjecture} holds. Then~$\ze(B_{E,P})$ is
bounded uniformly for~$\mathfrak F$, and the bound depends on~$c$
only. If, in addition to
Lang's Conjecture, either of the following conditions hold:
\begin{enumerate}
\item[\rm(1)] $P$ does not lie in
the (real) connected component of the identity;
\item[\rm(2)] $x(P)-x(Q)$
is a square,
\end{enumerate}
then~$\zo(B_{E,P})$ is bounded uniformly.
\end{theorem}

Infinite families satisfying the conditions of
Theorem~\ref{thetheorem2} are easy to manufacture.

\begin{example}\label{parametrizedexample} Fix~$T\in\mathbb N$,
$T>1$,
and let~$E$ denote the elliptic curve
$$
E:\quad y^2=x^3-T^2(T^2-1)x,
$$
together with the non-torsion point~$P=(1-T^2,1-T^2)$ and
the~$2$-torsion point~$Q=(0,0)$. Using the methods
in~{\rm\cite{MR2001i:11066}}, an explicit form of Lang's Conjecture is
provable for the family~$\mathfrak F=\{(E,P,Q)\}$. This gives an
example of case~{\rm(1)} in Theorem~\ref{thetheorem2}.
Taking~$P=(T^2,T^2)$ on the same curve yields an example of
case~{\rm(2)}.
\end{example}

\begin{example}\label{y2=x+1x-Tx-4T}
For all~$T>0$ consider the
curve
$$
E:\quad y^2=(x+1)(x-T)(x-4T),
$$
together with the non-torsion point~$P=(0,2T)$ and the~$2$-torsion
point~$Q=(-1,0)$. Lang's Conjecture holds for this family and, in
principle, the constant~$c$ can be computed explicitly. For this
family~{\rm(1)} in Theorem~\ref{thetheorem2} holds.
\end{example}

The proofs of the theorems
seem to need some form of Siegel's Theorem
on the finiteness of the number of integral points on the curve.
Indeed,~$\z(B_{E,P})$
being finite requires that~$B_n$ grows with~$n$.
There are effective versions of Siegel's Theorem,
however -- as far as we can see -- no routine
application of these will yield our results.
The strongest forms of
Siegel's Theorem are proved using elliptic transcendence theory.
These methods give good bounds in terms of the shape of error
terms and they work in
great generality. However, the dependence upon the
discriminant does not allow uniformity results -- also the size of
the constants gives excessively large estimates for the Zsigmondy
bound in particular cases. This is discussed
further after equation~\eqref{plugin1} below.

\subsection{Curves without rational~$2$-torsion}

The strongest results in the paper
require the presence of a rational~$2$-torsion
point. The following example illustrates how knowledge about the
odd Zsigmondy bound can outstrip that for the even bound when no
such point is present.

\begin{example}\label{y2+y=x3-x}
Consider the pair~$(E,P)$ with
$$
E:\quad y^2+y=x^3-x \mbox{ and }P=(0,0).
$$
The methods we describe allow a painless proof
that~$\zo(B_{E,P})=3$. Notice that in this case~$nP$ is integral
for~$n=1,2,3,4,6$ so we could not expect the bound to be any
smaller. However, we are unable to prove that the even Zsigmondy
bound is~$6$. Given any example where~$P$ does not lie in the real
connected component of the identity, the methods in this paper would
allow the odd Zsigmondy bound to be computed.
\end{example}

\begin{example}
The \emph{odd} terms of the sequence in Example~\ref{y2+y=x3-x}
comprise the Somos-4 sequence
$$
1,1,1,1,2,3,7,23,\dots.
$$
This sequence, which satisfies the bilinear recurrence
$$
u_nu_{n-4}=u_{n-1}u_{n-3}+u_{n-2}^2,
$$
was studied
by Somos~\cite{somoscrux}. By the bound for~$\zo$ in
Example~\ref{y2+y=x3-x}, every term of the Somos-4 sequence
beyond the fourth term has a primitive divisor.
\end{example}

For further results on Somos sequences, see the
papers~\cite{gale1},~\cite{MR93a:11012} and
the monograph~\cite[Sect.~1.1.17]{MR2004c:11015}.

Our final example is a family of curves for which knowledge about
the even Zsigmondy bound outstrips that for the odd bound. This is
included because it uses a new technique.

\begin{theorem}\label{y2=x3+T3+1}
Consider the pair~$(E,P)$ where
$$
E:\quad y^2=x^3+T^3+1 \mbox{ and }P=(-T,1).
$$
Then~$\ze(B_{E,P})$ is uniformly bounded
for all~$T>1$.
\end{theorem}

In the setting of Theorem~\ref{y2=x3+T3+1},
we are unable to prove such a
statement for the odd Zsigmondy bound.

The proof of Theorem~\ref{thetheorem}
is given in Section~\ref{sect:nonexplicitproofofthetheorem},
using a
sharpening of Silverman's original approach, together with
results of Bremner, Silverman and Tzanakis
concerning the difference between the na\"ive
height and the canonical height of a rational point on an elliptic
curve. In Section~\ref{explicitexamples} we will further
illustrate the method by explaining Examples~\ref{y2=x3-25x}
and~\ref{y2+y=x3-x}.
In Section~\ref{thetheorem2proof}, a proof
of Theorem~\ref{thetheorem2} will be given. Much of this is
routine and we will not labour it; however some explanation is
required for case~(1) in order to preserve the dependence of the
error term upon the discriminant. Theorem~\ref{y2=x3+T3+1} is proved in
Section~\ref{proofoftheoremy2=x3+t3+1}.

\section{Proof of Theorem~\ref{thetheorem}}
\label{sect:nonexplicitproofofthetheorem}

We begin with some basic facts about divisibility properties of
the sequence
$$
B_{E,P}=(B_n)_{n\ge1}.
$$

\begin{lemma}\label{ordthing}
Suppose~$p$
denotes any prime divisor of~$B_n$. Then
\begin{equation}\label{orderofmultiples}
{\ord}_p(B_{nk})={\ord}_p(B_n)+2{\ord}_p(k).
\end{equation}
\end{lemma}
This comes out of the development of
the~$p$-adic elliptic logarithm in~\cite{MR87g:11070}
and requires some local
analysis of elliptic curves. Note that the property of being a
divisibility sequence follows from~\eqref{orderofmultiples}. Indeed a
stronger property follows immediately.

\begin{lemma}\label{gcdthing}
For any~$m,n\in\mathbb N$
\begin{equation*}
\gcd(B_n,B_m)=B_{\gcd(m,n)}.
\end{equation*}
\end{lemma}

\begin{proof}
Let~$d=\gcd(m,n)$ and write~$m=kd$,~$n=\ell d$.
Then for any prime~$p$ dividing~$B_d$,
one of~$\ord_p(k)$ and~$\ord_p(\ell)$ must be zero.
By~\eqref{orderofmultiples},
$$
\ord_p(B_m)=\ord_p(B_d)+2\ord_p(k)\mbox{ and }
\ord_p(B_n)=\ord_p(B_d)+2\ord_p(\ell),
$$
so
\begin{eqnarray*}
\ord_p\left(\gcd(B_m,B_n)\right)&=&\min
\left\{\ord_p(B_d)+2\ord_p(k),\ord_p(B_d)+2\ord_p(\ell)\right\}\\
&=&\ord_p(B_d),
\end{eqnarray*}
so~$B_d\divides\gcd(B_n,B_m)$.
Conversely, if a prime~$p$ divides~$B_n$
and~$B_m$, then on the underlying elliptic
curve reduced modulo~$p$, $mP=nP=\mathcal{O}$, the
identity, hence~$dP=\mathcal{O}$
and so~$p\divides B_d$.
\end{proof}

These two lemmas will now be used to prove the fundamental
property shared by those terms~$B_n$ which do not have a primitive
divisor.

\begin{lemma}\label{fundprop}
If~$B_n$ does not have a primitive divisor
then
\begin{equation}\label{divpropinfundprop}
B_n\divides \prod_{p|n}p^2B_{{n}/{p}}.
\end{equation} If~\eqref{divpropinfundprop} holds, then
any primitive divisor of~$B_n$ divides~$n$.
\end{lemma}

\begin{proof}
Assume that~$B_n$ does not have a primitive divisor.
Let~$q$ be any prime, and~$p$ a prime dividing~$n$.
If~$\ord_q(B_{n/p})>0$ for some prime~$p\divides n$,
then by Lemma~\ref{ordthing}
$$
\ord_q(B_{n})=\ord_q(B_{n/p})+2\ord_q(p)
\le
\ord_q(B_{n/p})+2.
$$
If~$\ord_q(B_{n/p})=0$ for all primes~$p\divides n$
then~$q\notdivides B_n$. To see this,
notice that if~$q\divides B_n$ then by assumption~$q\divides B_m$
for some~$m\divides n$, hence~$q\divides B_{n/p}$ for some
prime~$p$, contradicting~$\ord_q(B_{n/p})=0$.

The partial converse follows in a similar way:
if~\eqref{divpropinfundprop} holds and~$q$ is a primitive divisor
of~$B_n$, then
$$
q\divides\prod_{p\divides n}p^2,
$$
so~$q\divides n$.
\end{proof}

Lemma~\ref{fundprop} will play a practical as well as a
theoretical role in the sequel. Our methods typically show
that~$\z(B_{E,P})\le C$ for some moderately large~$C$. The terms
with~$n\le C$ need to be checked to find the lowest bound. The
quadratic-exponential growth rate of the~$B_n$ means we wish to
avoid factorizing terms to do the checking. Lemma~\ref{fundprop}
is an easily implemented method for performing the check which is
factorization-free.

Finally, we gather some well-known facts about heights on elliptic curves.
Recall that~$P$ is a non-torsion point in~$E(\mathbb Q)$,
where the curve~$E$ is
\begin{equation*}
E:\quad y^2=x^3-T^2x,
\end{equation*}
with~$T\in \mathbb Z$
square-free.

Write~$h(\frac{a}{b})=\log\max\{\vert a\vert,\vert b\vert\}$
for the Weil height of
a rational number, so
$$
h(x(nP))=\log\max\{|A_n|,B_n\}.
$$

\begin{lemma}\label{heightstuff}
Let~$\hat h(P)$ denote the
global canonical height of~$P$. Then
\begin{equation}\label{maxbound}
n^2\hat h(P)-\textstyle\frac{1}{2}\log(T^2+1)-0.116\le h(x(nP))\le
n^2\hat h(P)+\log T+0.347,
\end{equation}
and
\begin{equation}\label{htlowerbound}
\hat h(P)\ge\textstyle\frac{1}{4}\log T.
\end{equation}
\end{lemma}

\begin{proof}
By~\cite[Eqn.~(15)]{MR2001i:11066},
for any point~$Q\in E(\mathbb Q)$,
$$
-0.347-\log T<\hat h(Q)-h(x(Q))<
\textstyle\frac{1}{2}\log(T^2+1)+
0.116
$$
(notice that the canonical height we are working with is
twice the value used in~\cite{MR2001i:11066}).
In particular,
\begin{eqnarray*}
h(x(nP))&\le&\hat h(nP)+\log T+0.347\\
&=&
n^2\hat h(P)+\log T+0.347
\end{eqnarray*}
and
\begin{eqnarray*}
h(x(nP))&\ge&\hat h(nP)-
\textstyle\frac{1}{2}\log(T^2+1)-0.116\\
&=&
n^2\hat h(P)-
\textstyle\frac{1}{2}\log(T^2+1)-0.116
\end{eqnarray*}
proving~\eqref{maxbound}.

The other result we call upon
also appeared in~\cite[Prop.~2.1]{MR2001i:11066}. If~$P$
denotes any non-torsion rational point on~$E$, then
\begin{equation*}
\textstyle\frac{1}{8}\log (2T^2)\leq
\hat h(P),
\end{equation*}
from which~\eqref{htlowerbound} is immediate.
\end{proof}

\proofof{Theorem~\ref{thetheorem}}
Assume that~$B_n$ does not have a primitive divisor.
Taking logarithms in Lemma~\ref{fundprop} gives
\begin{equation}\label{plugin1}
\log B_n \leq 2\sum_{p|n}\log p +\sum_{p|n}\log
B_{{n}/{p}}.
\end{equation}
The proof proceeds using various upper and lower estimates for~$\log B_k$
to make quantitative the observation that~\eqref{plugin1}
automatically bounds~$n$.

It is possible to use a deep general result from elliptic
transcendence theory to obtain a lower bound of
the form
\begin{equation}\label{sinnoudavidbound}
\log B_n\ge n^2\hat h(P)-\bigo(\log n\log \log n).
\end{equation}
Inserting this into~\eqref{plugin1} shows
that~$\z(B_{E,P})$ is finite because the
right-hand side is bounded by~$cn^2$ with~$c<1$.

Results of the form~\eqref{sinnoudavidbound}
have been obtained by
David~\cite{MR98f:11078}.
The form of the implied constant in~\eqref{sinnoudavidbound}
is given
explicitly in~\cite{MR95m:11056}. However, the shape of the constant is too
unwieldy for our purposes. For one thing, the dependence upon~$T$
comes as a power of~$\log T$ -- to obtain a uniformity result we
need it to be linear in~$\log T$. Another problem
is that the implied constants are
enormous. The quadratic-exponential
growth rate of the sequence~$B_{E,P}$ means that applying
this method would greatly complicate the computation of the
Zsigmondy bound.

Our approach is to use an inferior lower bound in
respect of the leading term: typically~$n^2\hat h(P)$ will be replaced by
three
quarters or even one quarter of this. However, the resulting error
term is more readily controlled.

By~\eqref{maxbound}, for any~$p\divides n$,
\begin{eqnarray}
\log B_{n/p}&\le&h(x(\textstyle\frac{n}{p}P))\nonumber\\
&\le&
\hat h(\textstyle\frac{n}{p}P)+\log T+0.347\nonumber\\
&=&\textstyle\frac{n^2}{p^2}\hat h(P)+\log T+0.347.\label{maxboundapplied}
\end{eqnarray}
We will call on three arithmetical
functions. Denote by~$\omega(n)$ the number of distinct prime
divisors of~$n$. Clearly
$$
\omega(n)\le\log n/\log 2 \le1.443\log n.
$$
Denote by~$\rho(n)$ the sum~$\sum_{p\vert n}\frac{1}{p^2}$
over prime divisors of~$n$. A calculation shows that
$$
\rho(n)\le0.453\mbox{ for all }n\ge1
$$
and, crucially,
$$
\rho(n)\le0.203\mbox{ for all odd }n\ge1.
$$
Finally, define
$$
\eta(n)=2\sum_{p\divides n}\log p.$$

Substituting~\eqref{maxboundapplied} into~\eqref{plugin1} gives
\begin{eqnarray}
\log B_n&\le&\eta(n)+\sum_{p\divides n}\left(
\textstyle\frac{n^2}{p^2}\hat h(P)+\log T+0.347\right)\nonumber\\
&\le&
\eta(n)+n^2\rho(n)\hat h(P)+
\omega(n)\left(\log T+0.347\right)\label{plugin2}
\end{eqnarray}

Assume first that~$n=2m$ is even. From
the duplication formula on the curve~$E$,
\begin{equation}\label{explicitdouble}
\frac{A_n}{B_n}=\frac{A_{2m}}{B_{2m}}=x(nP)
=x(2mP)=
\frac{(A_m^2+T^2B_m^2)^2}
{4A^{\vphantom{2}}_mB^{\vphantom{2}}_m(A_m^2-T^2B_m^2)}.
\end{equation}
It follows that
\begin{equation}\label{boundsb2mintermsofgcd}
B_{2m}=\frac{4A^{\vphantom{2}}_mB^{\vphantom{2}}_m(A_m^2-T^2B_m^2)}
{\gcd\left((A_m^2+T^2B_m^2)^2,
4A^{\vphantom{2}}_mB^{\vphantom{2}}_m(A_m^2-T^2B_m^2)\right)}.
\end{equation}
To bound the size of the greatest common divisor,
note that~$A_m$ and~$B_m$ are coprime by definition, and
recall that~$T$ is square-free.
We must allow for the possibility that~$4$ divides the numerator
in~\eqref{boundsb2mintermsofgcd}.
Now let~$p$ be an odd prime dividing the greatest common divisor.
Then
$$
p\divides A_m^2+T^2B_m^2
$$
and
$$
p\divides A_mB_m(A_m^2-T^2B_m^2),
$$
so~$p\divides A_m^3B_m$. Now~$p\divides B_m$ implies that~$p\divides A_m$, which
is impossible as~$A_m$ and~$B_m$ are coprime. So we deduce that~$p\divides A_m$
and hence~$p\divides T$.
Let
$$
\alpha=\ord_p(A_m^2+T^2B_m^2),
\quad\beta=\ord_p(A_m)\quad\mbox{ and }\quad\gamma=\ord_p(A_m^2-T^2B_m^2).
$$
If~$\beta\ge2$, then~$\alpha=\gamma=2$, so~$p$ divides the greatest
common divisor four times. If~$\beta=1$, then~$\gamma\ge2$ implies that~$\alpha=2$,
while~$\alpha\ge2$ implies that~$\gamma=2$. In all cases, it
follows that~$p$ divides the greatest common divisor no more than four
times. Thus
\begin{equation}\label{boundongcdindoubling}
{\gcd\left((A_m^2+T^2B_m^2)^2,
4A^{\vphantom{2}}_mB^{\vphantom{2}}_m(A_m^2-T^2B_m^2)\right)}
\le 4T^4.
\end{equation}

The greatest common divisor may also
be bounded using the following argument.
From~\eqref{maxbound}, trivial
estimates for the numerator and denominator
in~\eqref{explicitdouble} show that
the logarithm of each is bounded by
$
4\hat hm^2+\bigo(1),
$
with a uniform error. However~\eqref{maxbound} shows that~$\log
\max\{|A_{2m}|,B_{2m}\}$ is bounded below by~$4\hat hm^2-\bigo(\log T)$;
thus bounding the possible cancellation by a power of~$T$ as before.
For even~$n$ this approach is not needed, but we will make essential use of
it later for
one of the odd~$n$ cases.

From~\eqref{boundsb2mintermsofgcd} and~\eqref{boundongcdindoubling} we
deduce the important lower bound
$$
\frac{\vert
A^{\vphantom{2}}_mB^{\vphantom{2}}_m
(A_m^2-T^2B_m^2)\vert}{T^4}\le B_{2m},
$$
or in logarithmic form,
\begin{equation}
\log\vert A^{\vphantom{2}}_m\vert+\log B^{\vphantom{2}}_m+
\log\vert A_m^2-T^2B_m^2\vert
-4\log T\le\log B_{2m}=\log B_n.\label{punultimatebit?}
\end{equation}

\begin{lemma} For~$T\ge5$,
\begin{equation}
3\log \max \{|A_m|,B_m\}-\log T-0.693
\leq\log|A_m|+\log B_m+\log
|A_m^2-T^2B_m^2|.\label{finalbitofjigsaw}
\end{equation}
\end{lemma}

\begin{proof}
Let~$\alpha=\vert A_m\vert$ and~$\beta=\vert T\vert B_m$,
so that~\eqref{finalbitofjigsaw} follows from the
inequality
\begin{equation}\label{finalfinalbitofjigsaw}
\alpha\beta\vert\alpha-\beta\vert(\alpha+\beta)
\ge\frac{1}{2}\max\{\alpha,\beta\}^3.
\end{equation}
The expression in~\eqref{finalfinalbitofjigsaw} is symmetrical
in~$\alpha$ and~$\beta$, so assume without loss of generality
that~$\alpha>\beta$.

If~$\alpha\ge2\beta$ then
\begin{equation*}
\alpha\beta(\alpha-\beta)(\alpha+\beta)\ge
\alpha\cdot1\cdot\frac{1}{2}\alpha\cdot\alpha=\frac{1}{2}\alpha^3
\ge\frac{1}{2}\max\{\alpha,\beta\}^3.
\end{equation*}
If~$\beta<\alpha<2\beta$ then
\begin{equation*}
\alpha\beta(\alpha-\beta)(\alpha+\beta)\ge\alpha\cdot
\frac{1}{2}\alpha\cdot1\cdot\frac{3}{2}\alpha\ge
\frac{3}{4}\alpha^3\ge\frac{1}{2}\max\{\alpha,\beta\}^3.
\end{equation*}
\end{proof}

By~\eqref{punultimatebit?}
and~\eqref{finalbitofjigsaw},
\begin{eqnarray*}
\log B_n&\ge&\log\vert A_m\vert+\log B_m+\log\vert A_m^2-T^2B_m^2
\vert-4\log T\\
&\ge&
3\log\max\{\vert A_m\vert,B_m\}-5\log T-0.693\\
&=&
3h(x(mP))-5\log T-0.693,
\end{eqnarray*}
so by~\eqref{maxbound} and~\eqref{plugin2},
\begin{eqnarray}
\textstyle\frac{3}{4}n^2\hat h(P)-5\log T-\textstyle\frac{3}{2}\log(T^2+1)-
1.041&\le&
3h(x(mP))-5\log T-0.693
\nonumber\\
&\le&\log B_n\nonumber\\
&\le&
\eta(n)+n^2\rho(n)\hat h(P)\nonumber\\
&&\quad\quad\nonumber
+\omega(n)\left(\log T+0.347\right).
\end{eqnarray}
It follows that
\begin{eqnarray}
n^2\hat h(P)\left(\textstyle\frac{3}{4}-\rho(n)\right)
&\le&
\eta(n)+\omega(n)\left(\log T+0.347\right)\nonumber\\
&&\quad\quad+5\log T+
\textstyle\frac{3}{2}\log(T^2+1)+1.041.\label{boundintermsofht}
\end{eqnarray}
Recall that~$T\ge5$ so~$\log T>1.609$ and hence
\begin{equation}\label{forwardrefer}
\frac{\log(T^2+1)}{\log T}\le
2.0244.
\end{equation}
Apply~\eqref{htlowerbound}
to~\eqref{boundintermsofht}, divide
through by~$\log T$,
and apply~\eqref{forwardrefer}
to deduce that
$$
n^2\left(\vphantom{\textstyle\frac{3}{4}}
\textstyle\frac{3}{4}-\rho(n)\right)\le
4\left(\vphantom{A^A}0.621\eta(n)+1.216\omega(n)
+9.0776
\right).
$$
This implies that~$n\le11$, so~$\ze(B_{E,P})\le10$.

The bound obtained so far (when~$n$ is even) takes a similar form
in general. Assume that~$x(P)<0$ and~$n$ is odd. If~$B_n\ge |A_n|$
then
\begin{equation}\label{lastchangemaybe!}
\log B_n \ge h(x(nP))\ge n^2\hat
h(P)-\textstyle\frac{1}{2}\log(T^2+1)-0.116.
\end{equation}
If~$B_n<|A_n|$,
use the fact that if~$n$ is odd then~$x(nP)<0$,
therefore
$$
-T\leq x(nP)< 0.
$$
Thus~$|A_n/B_n|\le T$, so
$$
\log |A_n|-\log T\le \log B_n.
$$
Therefore
\begin{eqnarray*}
\log B_n&\ge&h(x(nP))-\log T\\
&\ge&n^2\hat h(P)-\textstyle\frac{1}{2}\log(T^2+1)-\log T-0.116.
\end{eqnarray*}
This lower bound, being
smaller than the one in~\eqref{lastchangemaybe!},
covers both cases.
By~\eqref{maxbound} and~\eqref{plugin1},
\begin{eqnarray*}
n^2\hat h(P)-\textstyle\frac{1}{2}\log(T^2+1)-\log T-0.116&\le&
\log B_n\\
&\le&
\eta(n)+n^2\rho(n)\hat h(P)\\
&&\medspace\medspace\medspace\medspace+\omega(n)\left(\log T
+0.347\right).
\end{eqnarray*}
The bound~\eqref{htlowerbound} then
implies that
$$
n^2(1-\rho(n))\le4\left(\vphantom{A^A}
0.621\eta(n)+
1.216\omega(n)+2.085\right),
$$
using~\eqref{forwardrefer} again.

It follows (for odd~$n$) that~$n\le3$, so~$\zo(B_{E,P})\le3$. This dramatic
improvement in the size of the bound is mainly accounted for by
the fact that~$\rho(n)\le0.203$ for all odd~$n$, and the very good
lower bound for $\log B_n$. The fact that the bound for~$\rho(n)$
over odd~$n$ is strictly smaller than~$\frac{1}{4}$ will play a
critical role later.

Finally, assume that~$x(P)$ is a square.
For this
part of Theorem~\ref{thetheorem}, we are going to use the fact
that~$x(nP)$ is a square for all~$n\in \mathbb N$. This follows
from the proof of the Weak Mordell Theorem: the map~$E(\mathbb
Q)\rightarrow \mathbb Q^*/\mathbb Q^{*2}$ given by$$P\mapsto
x(P)\mathbb Q^{*2} \mbox{ and } (0,0)\mapsto -\mathbb Q^{*2}$$ is
a group homomorphism. Write
\begin{equation*}
nP=\left(\frac{A_n}{B_n^{\vphantom{3/2}}},
\frac{C_n}{B_n^{3/2}}\right)\negmedspace.
\end{equation*}
Assume that~$n=2m+1$ is odd and write
$$
nP=mP+(m+1)P.
$$
Then
$$
x(nP)=\frac{A_{2m+1}}{B_{2m+1}}=
\left(\frac{y((m+1)P)+y(mP)}{x((m+1)P)-x(mP)}
\right)^2-x(mP)-x((m+1)P).$$
Inserting the explicit form of~$nP=mP+(m+1)P$
into this
formula yields
\begin{equation}\label{oddsecondcase}
\frac{(A_mA_{m+1}-T_{\vphantom{1}}^2B_mB_{m+1})
(A_mB_{m+1}+A_{m+1}B_m)-2C_mC_{m+1}B_{m\vphantom{1}}^{1/2}B^{1/2}_{m+1}}{(A_{m+1}B_m
-A_mB_{m+1})^2}
\end{equation}
for~$x(nP)$.
Once again we wish to bound the possible size of the greatest
common divisor of the
numerator~$N$ and the denominator~$D$ in~\eqref{oddsecondcase}.
An additional complication here is the appearance of terms
arising from~$y(P)$. Since~$nP$ lies on the curve~$y^2=x^3-T^2x$,
$$
C_n^2=A_n^3-T^2A_nB_n^2.
$$
It follows that
\begin{eqnarray}
\log\vert C_n\vert&\le&\textstyle
\frac{1}{2}\left(\log 2+
\max\log\{|A_n|^3,T^2|A_n|^{\vphantom{2}}B_n^2\}\right)\nonumber\\
&\le&\textstyle\frac{1}{2}
\left(\log2+3n^2\hat h(P)+5\log T+1.041\right)\negmedspace.\label{boundonc}
\end{eqnarray}
Now write
$$
\alpha=(A_mA_{m+1}-T^2B_mB_{m+1})
(A_mB_{m+1}+A_{m+1}B_m)
$$
and
$$
\beta=2C^{\vphantom{1/2}}_mC^{\vphantom{1/2}}_{m+1}B_{m\vphantom{1}}^{1/2}B^{1/2}_{m+1}.
$$
By using~\eqref{maxbound} and~\eqref{boundonc},
\begin{eqnarray}
\log\vert\alpha\vert&\le&\log4+\log\max\{|A_mA_{m+1}|,T^2B_mB_{m+1}\}\nonumber\\
&&\quad\quad\quad\quad+
\log\max\{|A_m|B_{m+1},|A_{m+1}|B_m\}\nonumber\\
&\le&
(4m^2+4m+2)\hat h(P)+
6\log T+2.775\nonumber
\end{eqnarray}
and
\begin{equation*}
\log\vert\beta\vert\le(4m^2+4m+2)\hat h(P)+6\log T+2.775.
\end{equation*}
Thus the numerator and denominator
of~\eqref{oddsecondcase} satisfy
\begin{eqnarray}
\max\{\log\vert N\vert,\log\vert D\vert\}&\le&
\log2+\log\max\{\vert\alpha\vert,\vert\beta\vert\}\nonumber\\
&\le&
(4m^2+4m+2)\hat h(P)+6\log T+3.469.
\end{eqnarray}
On the other hand, by the lower bound in~\eqref{maxbound},
\begin{eqnarray}
\max\{\log\vert A_n\vert,\log B_n\}
&\ge&
n^2\hat h(P)-\textstyle\frac{1}{2}\log(T^2+1)-0.116\nonumber\\
&=&
(4m^2+4m+1)\hat h(P)-\textstyle\frac{1}{2}\log(T^2+1)-0.116.\nonumber
\end{eqnarray}
It follows that
$$
\gcd(N,D)\le
\hat h(P)+6\log T+\textstyle\frac{1}{2}\log(T^2+1)+3.584,
$$
so by~\eqref{oddsecondcase} and~\eqref{plugin2}
\begin{eqnarray}
2\log\left(A_{m+1}B_m-\right.&&\left.
\negmedspace\negmedspace\negmedspace\negmedspace A_mB_{m+1}\right)-
\hat h(P)-6\log T-
\textstyle\frac{1}{2}\log(T^2+1)-3.584\nonumber\\
&<&\log B_n\nonumber\\
&<&
\eta(n)+n^2\rho(n)\hat h(P)+
\omega(n)\left(\log T+0.347\right).\label{plugin3oddcase}
\end{eqnarray}
Now by assumption~$A_m,A_{m+1},B_m$ and~$B_{m+1}$ are all
squares; write~$A_{*}=a_{*}^2$ and~$B_{*}=b_{*}^2$ with~$a_{*},b_{*}>0$.
Then
\begin{eqnarray}
\max\{\log\vert a_{m+1}\vert,\vert b_{m+1}\vert\}&\le&
\log\left(\vert a_{m+1}\vert+\vert b_{m+1}\vert\right)\nonumber\\
&\le&\log\left(\vert a_{m+1}b_m\vert+\vert a_mb_{m+1}\vert\right)\nonumber\\
&\le&\log\left\vert a_{m+1}^2b_m^2-a_m^2b_{m+1}^2\right\vert,
\end{eqnarray}
so by~\eqref{plugin3oddcase}
\begin{eqnarray}
h\left((m+1)P\right)&=&\max\{\log A_{m+1},B_{m+1}\}\nonumber\\
&\le&
\eta(n)+(n^2+1)\rho(n)\hat h(P)+\omega(n)\left(
\log T+0.347\right)\nonumber\\
&&\quad\quad\quad\quad+6\log T+\textstyle\frac{1}{2}\log(T^2+1)
+3.584.\nonumber
\end{eqnarray}
Using~\eqref{maxbound},~\eqref{htlowerbound} and the
assumption that~$T\ge5$, this shows that
\begin{equation}\label{looksbad}
\textstyle\frac{1}{4}(n+1)^2-(n^2+1)\rho(n)\le
4\left(\vphantom{A^A}
0.621\eta(n)+10.596+1.216\omega(n)
\right).
\end{equation}
It is not clear that the left-hand side
of~\eqref{looksbad} grows at all. However, as noted earlier, for
odd~$n$ we have~$\rho(n)<0.203<\frac{1}{4}$, so the left-hand side
of~\eqref{looksbad} grows at least like~$0.047n^2$ for odd~$n$.
Thus~(\ref{looksbad}) does bound~$n$. Indeed~\eqref{looksbad}
implies that~$n\le21$, showing that~$\zo(B_{E,P})\le21$.
\endproofof

\section{Explicit Examples}\label{explicitexamples}
Theorem \ref{thetheorem} supplies such good bounds that the
remaining cases can be checked using Lemma~\ref{fundprop}.
Inserting explicit values for the canonical heights in specific
examples reduces the checking even further. From the proof in
Section~\ref{sect:nonexplicitproofofthetheorem} we have the
following inequalities under the assumption that~$B_n$ does not
have a primitive divisor. If~$x(P)<0$ and~$n$ is odd, then
\begin{eqnarray}
n^2\hat h(P)\left(1-\rho(n)\right)&\le&\eta(n)+\omega(n)\log T+0.347\omega(n)\nonumber\\
&&\quad\quad\quad\quad+\textstyle\frac{1}{2}\log(T^2+1)
+\log T+0.116;\label{extract1}
\end{eqnarray}
whilst if~$n$ is even, then
\begin{eqnarray}
n^2\hat h(P)\left(\textstyle\frac{3}{4}-\rho(n)\right)
&\le&
\eta(n)+\omega(n)\left(\log T+0.347\right)\nonumber\\
&&\quad\quad+5\log T+
\textstyle\frac{3}{2}
\log(T^2+1)+1.041.\label{extract2}
\end{eqnarray}

{\sc Example~\ref{y2=x3-25x}.}
Here~$T=5$ and the canonical height of~$P=(-4,6)$ is given by $
\hat h(P)=1.899\dots. $ Theorem~\ref{thetheorem}
predicts~$\ze(B_{E,P})\le 12$. Using
Lemma~\ref{fundprop}, the checking of the remaining cases is
quick. Inserting the explicit estimate for $\hat h(P)$ reduces
this calculation still further. Assuming that~$B_n$ does not have
a primitive divisor, ~\eqref{extract1} and~\eqref{extract2} imply
$\zo(B_{E,P})=1$ and~$\ze(B_{E,P})\le8.$ The remaining cases can
easily be checked almost by hand, but certainly using
Lemma~\ref{fundprop}.

\smallskip

{\sc Example~\ref{y2+y=x3-x}.}
This is proved in similar fashion to Example~\ref{y2=x3-25x} so it
is not discussed it in detail.

\section{Proof of Theorem~\ref{thetheorem2}}\label{thetheorem2proof}
Suppose without loss of
generality that~$Q=(0,0)$ in every case, since translation
preserves both the discriminant of the curve and the kind of
result sought. Assume the defining equation for~$E$ has the form
$$
E:\quad y^2=x(x^2+ax+b)=x(x-r_1)(x-r_2).
$$
The discriminant~$\Delta=\Delta(E)$ of the
curve is given by
\begin{equation}\label{equationfordiscinthetheorem2}
\Delta=(r_1r_2(r_1-r_2))^2.
\end{equation}

\begin{lemma}
\begin{equation}\label{rootsanddiscbound}
\max\{|\log |r_1||,|\log |r_2||\}\le\frac{3}{2}\log|\Delta|.
\end{equation}
\end{lemma}

\begin{proof}
Without loss of generality we may assume that~$\vert r_1\vert
\le\vert r_2\vert$. It follows that~$\vert r_2\vert\ge1$
and~$\vert r_1\vert\ge\frac{1}{\vert r_2\vert}$, so
\begin{equation*}
\max\{\vert\log\vert r_1\vert\vert,\vert\log\vert r_2\vert\vert\}=
\log\vert r_2\vert.
\end{equation*}
If~$\vert r_1\vert\le\frac{1}{2}\vert r_2\vert$ then
\begin{equation*}
\Delta=\vert r_2r_2\left(\frac{r_1}{r_2}-1\right)\vert^2
\ge
\frac{b^2\vert r_2^2\vert}{4}\ge
\frac{r_2^2\vert}{4}
\end{equation*}
so~$\vert r_2\vert\le 2\sqrt{\vert\Delta\vert}.$

Assume now
that~$\vert r_1\vert>\frac{1}{2}\vert r_2\vert$.
Now
\begin{equation*}
\vert r_1-r_2\vert=\sqrt{\vert a^2-4b\vert}\ge1,
\end{equation*}
so
\begin{equation*}
\vert\Delta\vert\ge r_1^2r_2^2\ge\frac{1}{4}\vert r_2\vert^4,
\end{equation*}
and thus~$\vert r_2\vert\le\left(4\vert\Delta\vert\right)^{1/4}.$
Since~$\vert\Delta\vert\ge3$, this completes the proof.
\end{proof}

In the situation of Theorem~\ref{thetheorem2}, we need
a bound of the form
\begin{equation}\label{clarifythis}
|\hat h(P)-h(P)|\leq c\log \Delta,
\end{equation}
and this follows from the result in~\cite{MR91d:11063}
which bounds~$|\hat h(P)-h(P)|$ in terms of the height of
the~$j$-invariant
(and hence the height of the discriminant) of the curve.

\proofof{Theorem~\ref{thetheorem2}}
In the even case, writing~$n=2m$ and applying the duplication
formula shows that
$$\log |A_m|+\log B_m + \log |A_m^2+aA_mB_m+bB_m^2|-\bigo(\log \Delta)\leq \log
B_n.
$$
If~$|A_m^2+aA_mB_m+bB_m^2|\geq |A_mB_m|$ then
\begin{eqnarray*}
\log B_n&\ge&2h(mP)-\bigo(\log\Delta)\\
&\ge&\textstyle\frac{1}{2}\hat{h}(P)-\bigo(\log\Delta)
\end{eqnarray*}
by~\eqref{clarifythis}.
On the other hand, using the same argument as before shows
$$
\log |A_m|-\log B_m = \bigo(\max \{|\log |r_1||,|\log |r_2||\}=
\bigo(\log \Delta)
$$
by~\eqref{rootsanddiscbound}. This gives an analog of
the inequality~\eqref{boundintermsofht}, and the proof
proceeds as before.

In case~(2), the argument for the odd Zsigmondy bound is essentially
identical to that given before. In case~(1) the existence of two
connected components requires there to be three real~$2$-torsion
points; there are
then various cases to consider depending upon the
signs and relative sizes of the roots, and these can be summarized
as follows. Notice
first that
$$
\log\vert A_n/B_n\vert\leq \max \{|\log |r_1||,|\log |r_2||,\log
|r_1-r_2|\}.
$$
Each of the terms on the right is~$\bigo(\log \Delta)$ and
$$
hn^2-\bigo(\log \Delta)\leq \log B_n.
$$
The proof is completed exactly as before.
\endproofof

\section{Proof of Theorem~\ref{y2=x3+T3+1}}\label{proofoftheoremy2=x3+t3+1}
This may be shown using strong results of
Bennett~\cite{MR98c:11070},~\cite{MR98d:11082}
on Diophantine approximation in addition to the methods
of Section~\ref{sect:nonexplicitproofofthetheorem}.
Writing~$n=2m$ as
usual, the crucial point is to find an explicit estimate for
$$
B^{\vphantom{3}}_m\left|A_m^3+(T^3+1)B_m^3\right|.
$$
If~$A_m/B_m$ is bounded away from ~$\theta=(T^3+1)^{\frac{1}{3}}$
then we can proceed as before without difficulty. Otherwise, we
need some kind of explicit lower bound from Diophantine
approximation, of the form
$$
\frac{a}{q^{\lambda}}<\left|\theta-\frac{p}{q}\right|
$$
for all rationals~$p/q$ in lowest terms. Probably the best results
of this kind have been found by Bennett~\cite{MR98c:11070},~\cite{MR98d:11082}.
Applying these
estimates shows we may take
$$
\log a=\bigo(\log T)
$$
and
\begin{equation}\label{bennetfunction}
\lambda=1+\frac{2\log (\sqrt {T^3} + \sqrt
{T^3+1})+\log(3\sqrt 3/2)}{2\log (\sqrt {T^3} + \sqrt
{T^3+1})-\log(3\sqrt 3/2)},
\end{equation}
where all implied constants are explicit and uniform.
The right-hand side of~\eqref{bennetfunction}
is decreasing in~$T$ and converges to~$2$ as~$T\to\infty$.
For the methods used here,
we need~$\lambda<2.188$ and for this~$T$ needs to be at least~$26$.
Inserting this data into our machine yields an
inequality of the form
$$
\hat h\left(\vphantom{A^A}
0.047+\bigo(1/\log T)\right)n^2<2\log n + \bigo(\log T).
$$
Finally, the
canonical height of~$P$ satisfies
$$
\hat h=\hat h(P)\sim\textstyle\frac{1}{2}\log T.
$$
Using the same methods as in~\cite{MR2001i:11066}, it is possible to give an
explicit, positive lower bound for~$\hat h(P)/\log T$ and the
uniformity result follows. For this class of examples we have not
tried to state the most explicit result possible.


\end{document}